\documentclass[a4paper,12pt]{amsart}
\usepackage{amssymb}
\usepackage{eucal}
                                                                            
\usepackage{amsmath}
\usepackage{amscd}
\usepackage{amsthm}
                                                                                
\usepackage[T1]{fontenc}

\usepackage[latin1]{inputenc}

\theoremstyle{plain}
\newtheorem{thm}{Theorem}[section]

\newtheorem{lem}[thm]{Lemma}

\theoremstyle{definition}

\theoremstyle{remark}

\newcommand{\cal}{\mathcal}

\newcommand{\GL}{\mathrm{GL}}
\newcommand{\goth}{\mathfrak}

\newcommand{\id}{\mathrm{Id}}

\newcommand{\M}{\mathrm{M}}

\newcommand{\Quot}{\mathrm{Quot}}

\newcommand{\Rep}{\mathrm{Rep}}

\newcommand{\Spec}{\mathrm{Spec}}

\newcommand{\al}{\alpha}

\newcommand{\De}{\Delta}

\newcommand{\ph}{\varphi}
\newcommand{\si}{\sigma}

\newcommand{\la}{\lambda}

\begin{document}
\title[Differential Galois group.]{A note about rational representations of differential Galois groups.}
\author{Marc Reversat}

\address{Intitut de Math\'ematiques de Toulouse, Universit\'e Paul Sabatier, 118 route de Narbonne, 31062 Toulouse cedex 9, France.}

\email{marc.reversat@math.ups-tlse.fr}
\date {\today }


\begin{abstract}
We give a description of the rational representations of the differential Galois group of a Picard-Vessiot extension.

\end{abstract}


\maketitle

\subjclass{}
\keywords{}


\thanks{}
\dedicatory{}

\tableofcontents


We give a description of the rational representations of the differential Galois group of a Picard-Vessiot extension 
(theorems \ref{main thm} and \ref{2main thm}). This gives a new description of the differential Galois correspondence.
More results will be given for abelian differential extensions in a forthcoming paper \cite{vdP-R}, especially the analog of the Artin correspondence.

\section{On representations of Galois groups of Picard-Vessiot extensions.}
In all this note we consider differential fields with algebraically closed fields of constants, denoted by $C$.
The derivative will be denoted by a dash.\par

Let $K$ be a differential field. Let $n\geq 1$ be an integer, $\M_n(K)$ and $\GL_n(K)$ are the usual notations for algebra and group of $n\times n$ matrices with entries in $K$. 
The group $\GL_n(K)$ acts on $\M_n(K)$ by the following rule:
	\[
\begin{array}{ccc}
\GL_n(K)\times \M_n(K) & \longrightarrow & 	\M_n(K)\\
\left(U, A\right) & \longmapsto & U'U^{-1}+UAU^{-1}
\end{array}
\]
where if $U=(u_{i,j})_{1\leq i,j\leq n}$, then $U'=(u'_{i,j})_{1\leq i,j\leq n}$. This action can be defined in an other way, maybe more comprehensible. Consider the group
	\[H_n(K):=\M_n(K)\times \GL_n(K),
\]
the law being defined by the following formula: for all $A,B$ in $\M_n(K)$ and all $F,G\in \GL_n(K)$
	\[(A,F)(B,G)=(A+FBF^{-1},FG).
\]
It admits the subgroups 
	\[\De_n(K):=\left\{(U'U^{-1},U) \ / \ U\in \GL_n(K)\right\},
\]
$\left\{0\right\}\times \GL_n(K)$ and $\M_n(K)\times \left\{1\right\}$, this last one being normal. We set
	\[Z_n(K):=\De_n(K)\backslash H_n(K)/\left(\left\{0\right\}\times \GL_n(K)\right).
\] 
The group $\GL_n(K)$ acts on $\M_n(K)\times \left\{1\right\}$ by the following rule:
	\[
\begin{array}{ccc}
\GL_n(K)\times \left(\M_n(K)\times \left\{1\right\}\right) & \longrightarrow & 	\left(\M_n(K)\times                                                                                                       \left\{1\right\}\right)\\
\left(U, (A, 1)\right) & \longmapsto & (U'U^{-1}, U)(A, 1)(0,U^{-1})\\
{} & {} & =(U'U^{-1}+UAU^{-1}, 1)
\end{array}
\]
With the identification $\M_n(K)=\left(\M_n(K)\times \left\{1\right\}\right)$ and
inclusion \newline
$\left(\M_n(K)\times \left\{1\right\}\right)\subset H_n(K)$, it induces a canonical bijection
\begin{equation} \label{def Zn}
	\GL_n(K)\backslash \M_n(K)\simeq 
\De_n(K)\backslash H_n(K)/\left(\left\{0\right\}\times \GL_n(K)\right)=	Z_n(K).
\end{equation}
We will use below these two definitions of $Z_n(K)$.\par
 
 Let $A$ be in $\M_n(K)$ and $F$ be in $\GL_n(K)$, we denote by $[(A,F)]$, resp. $[A]$, the class in $Z_n(K)$ of $(A,F)\in H_n(K)$, resp. of $A\in \M_n(K)$. \par
 
Let $L/K$ be a Picard-Vessiot extension. The inclusion $H_n(K)\subseteq H_n(L)$ gives rise to a map
$\al(L/K):Z_n(K)\rightarrow Z_n(L)$. We set
	\[Z_n(L/K):=\left\{a\in Z_n(K) \ / \ \al(L/K)(a)=[0]\right\}.
\] \par

For any group $G$ we denote by $\Rep_n(G)$ the set of equivalent classes of representations of $G$ in $\GL_n(C)$,
if $G=dGal(L/K)$ is the differential group of $L$ over $K$, we set $\Rep_n(L/K):=\Rep_n(G)$.
 
\begin{thm} \label{main thm}
Let $L/K$ be a Picard-Vessiot extension, then there exists a natural bijection between $Z_n(L/K)$ and $\Rep_n(L/K)$.
\end{thm}

\begin{proof} First of all we recall some facts that are of main importance in our proof. \par

\subsubsection{The representation $c_A$.} \label{c_A}
Consider the differential equation $Y'=AY$, with $A\in \M_m(K)$, and let
$E/K$ be a corresponding Picard-Vessiot extension, i.e. $E$ is generated over $K$ by the coefficients of a fundamental matrix $F_A$ of the equation. The rational representation $c_A$ is
	\[
\begin{array}{ccc}
	dGal(E/K) & \longrightarrow & \GL_m(C) \\
	\si       & \longmapsto     & c_A(\si )
\end{array}
\]
where $c_A(\si )$ is such that $\si (F_A)=F_Ac_A(\si )$. Note that $c_A$ depends only on the class $\left[A\right]$ of
$(A,1)$ in $Z_m(K)$, because if $B=U'U^{-1}+UAU^{-1}$ with $U\in \GL_m(K)$, a fundamental matrix of the equation
$Y'=BY$ is $UF_A$ and we see that $c_B=c_A$. Note also that an other fundamemtal matrix is of the form $F_A\gamma $, with $\gamma \in \GL_n(C)$, then it gives the representation
$\gamma^{-1}c_A \gamma$ equivalent to $c_A$.\par

We will write equivalently $c_{A}$, $c_{[A]}$ or $c_E$ for this class of representations.

\subsubsection{The Galois group $dGal(E/K)$.} \label{Galois group}
Let $A\in \M_m(K)$ and 
\[
	R=K[(X_{i,j})_{1\leq i,j\leq m}, (\det )^{-1}]/q=K[(x_{i,j})_{1\leq i,j\leq m}]
\]
be a Picard-Vessiot ring over $K$ for the equation $Y'=AY$. In these formulas the $X_{i,j}$ are indeterminates, the ring $K\left[(X_{i,j})_{1\leq i,j\leq m} \right]$ is equipped by the derivation satisfying 
$(X'_{i,j})_{1\leq i,j\leq m}=A(X_{i,j})_{1\leq i,j\leq m}$, ``$\det $'' is the determinant of the matrix $(X_{i,j})_{1\leq i,j\leq m}$, $q$ is a maximal differential ideal and $x_{i,j}$ is the image of $X_{i,j}$. 
Let $E=\Quot (R)$ and ${\goth U}=dGal(E/K)$, consider 
\[
	K\left[(X_{i,j})_{1\leq i,j\leq m}, (\det)^{-1} \right]\subseteq 
	E\left[(X_{i,j})_{1\leq i,j\leq m}, (\det)^{-1} \right]
\]
\[	=E\left[(Y_{i,j})_{1\leq i,j\leq m}, (\det)^{-1} \right]\supseteq 
	C\left[(Y_{i,j})_{1\leq i,j\leq m}, (\det)^{-1} \right],
\]
where $(Y_{i,j})_{1\leq i,j\leq m}$ is defined by 
$(X_{i,j})_{1\leq i,j\leq m}=(x_{i,j})_{1\leq i,j\leq m}(Y_{i,j})_{1\leq i,j\leq m}$.
Note that $Y'_{i,j}=0$. We know that
\[
	{\goth U}=dGal(E/K)=\Spec C\left[(Y_{i,j})_{1\leq i,j\leq m}, (\det)^{-1} \right]/J
\]
where $J=qE\left[(Y_{i,j})_{1\leq i,j\leq m}, (\det)^{-1} \right]\cap 
C\left[(Y_{i,j})_{1\leq i,j\leq m}, (\det)^{-1} \right]$ (\cite{vdP-S} proof of prop. 1.24 or the beginning of
\S 1.5).
We denote by $y_{i,j}$ the image of $Y_{i,j}$, then we have
	\[{\goth U}=dGal(E/K)=\Spec C\left[(y_{i,j})_{1\leq i,j\leq m} \right]
\]\par

\subsubsection{${\goth U}=dGal(E/K)$ as a torsor.} \label{torsor}
We continue with the previous notations. Set ${\cal T}=\Spec R$, we know that 
${\cal T}$ is an ${\goth U}$-torsor over $K$ (\cite{vdP-S} theorem 1.30), moreover, we know that there exists a finite extension $\widetilde{K}$ of $K$ such that 
${\cal T}\times_K\widetilde{K}=\Spec\left(R\otimes_K\widetilde{K}\right)$ is a trivial ${\goth U}$-torsor over 
$\widetilde{K}$ (\cite{vdP-S} cor. 1.31), this means that there exists 
$\underline{b}\in {\cal T}(\widetilde{K})$ such that the following map is an isomorphism of $\widetilde{K}$-schemes
\[
\begin{array}{rccl}
\psi:  &	{\goth U}\times_{C}\widetilde{K} & \longrightarrow & {\cal T}\times_{K}\widetilde{K}  \\
{} &	(c_{i,j})_{1\leq i,j\leq m} & \longmapsto & \underline{b}(c_{i,j})_{1\leq i,j\leq m}
\end{array}
\]
($\underline{b}$ can be seen as a matrix, on the right this is a product of matrices; see the definition of $R$ above). \par

\subsubsection{Galois actions.} \label{Galois action}
Let $\si $ be an element of ${\goth U}=dGal(E/K)$, the action of $\si $ on $R$ is given by the images of the $x_{i,j}$, which are defined by the matrix formula $(\si (x_{i,j}))=(x_{i,j})c_E(\si )$. We denote by $\si^{\flat}$ the morphism induces by $\si $ on ${\cal T}$ or on ${\cal T}\times_{K}\widetilde{K}$, this is the action of 
${\goth U}$ which defines the torsor  structure.
An element of ${\cal T}(\widetilde{K})$ can be represented by a matrix 
$\underline{a}=(a_{i,j})_{1\leq i,j\leq m}$ with $a_{i,j}$ in $\widetilde{K}$, its image is 
$\si^{\flat}(\underline{a})=\underline{a}c_E(\si )$. For any $\si $ in ${\goth U}$ denote by $\la_{\si}$ the right translation on ${\goth U}$ by $\si $, i.e. 
	\[
\begin{array}{rccl}
\la_{\si}: & {\goth U} & \rightarrow & {\goth U} \\
{}         & \tau      & \mapsto     & \tau \si 
\end{array}
\]
Write again $\la_{\si}$ for 
$\la_{\si}\times \id_{\widetilde{K}} :{\goth U}\times_C \widetilde{K} \rightarrow {\goth U}\times_C\widetilde{K} $,
then the morphism $\psi $ of (\ref{torsor}) is equivariant, this means that for any $\si \in {\goth U}$ the following diagram is commutative
\[
\begin{array}{ccc}
{\goth U}\times_C\widetilde{K} & \stackrel{\psi}{\rightarrow}	& {\cal T}\times_K\widetilde{K} \\
\la_{\si }\downarrow & {} & \downarrow \si^{\flat} \\
{\goth U}\times_C\widetilde{K} & \stackrel{\psi}{\rightarrow}	& {\cal T}\times_K\widetilde{K} 
\end{array}
\]
\par

\emph{Proof of the theorem, the map $Z_n(L/K)\rightarrow R_n(L/K)$.}  \par

Let $A\in \M_n(K)$ such that $[A]\in Z_n(L/K)$, then there exists $U\in \GL_n(L)$ such that
	\[A=U'U^{-1},
\]
this means that $U$ is a fundamental matrix of the equation $Y'=AY$, as it is with entries in $L$, it exists a differential subextension $E$ of $L$ which is a Picard-Vessiot extension for the equation $Y'=AY$. Denote by  
$\rho_A $ the representation
\begin{equation} \label{2}
	\rho_A :dGal(L/K)\stackrel{{\mathrm{restriction}}}{\longrightarrow} dGal(E/K)
		\stackrel{c_A}{\longrightarrow} \GL_n(C).
\end{equation} \par

Now we prove that this representation $\rho_A $ does not depend on the class of $A$ in $Z_n(K)$ and of the choice of
$U\in \GL_n(L)$ such that  $A=U'U^{-1}$.\par

Let $B\in \M_n(K)$ such that $[B]=[A]$ in $Z_n(K)$, then there exists $W,T\in \GL_n(K)$ such that
	\[(B,1)=(W'W^{-1},W)(A,1)(0,T),
\]
it follows that $B=W'W^{-1}+WAW^{-1}$, this means that $WU$ is a fundamental matrix of the equation
$Y'=BY$. We see that $\rho_A =\rho_B$, and we denote this representation by $\rho_{[A]} $. \par

Let $V\in \GL_n(L)$ such that $A=U'U^{-1}=V'V^{-1}$, then we see that $\left(V^{-1}U\right)'=0$, this means that
there exists $\gamma \in \GL_n()$ such that $U=V\gamma$ and the two representations define as before in (\ref{2}) are conjugate.\par

Then to each element $[A]$ of $Z_n(L/K)$ we have associated the element $\rho_{[A]} $ of $\Rep_n(L/K)$. \par

\emph{Proof of the theorem, the map $R_n(L/K)\rightarrow Z_n(L/K)$.}  \par

Let $\rho :dGal(L/K)\rightarrow \GL_n(C)$ be a rational representation. Let $E$ be the fixed field of $\ker \rho$, we set ${\goth U}=dGal(E/K)$ and we denote again by $\rho$ the representation ${\goth U}\hookrightarrow \GL_n(C)$ coming from the given one. The field $E$ is a Picard-Vessiot extension corresponding to an equation $Y'=AY$, with 
$A\in \M_m(K)$. Our aim is to prove that one can chose $A$ in $\M_n(K)$, i.e. $m=n$, and that this gives the inverse map of $[A]\mapsto \rho_{[A]}$.\par 

We use the previous notations and descriptions of $E$, $R$, ${\goth U}$, ${\cal T}$ etc. We set
$\GL_n(C)=\Spec C[(T_{r,s})_{1\leq r,s\leq n}, (\det )^{-1}]$, let
\[
\rho^{\sharp}:C\left[(T_{r,s})_{1\leq r,s\leq n}, (\det)^{-1} \right]\longrightarrow
	C\left[(Y_{i,j})_{1\leq i,j\leq m}, (\det)^{-1} \right]/J
\]
be the comorphism of $\rho :{\goth U}\hookrightarrow \GL_n(C)$; $\rho^{\sharp}$ is onto. Set $I=\ker (\rho^{\sharp})$,
then we have an isomorphism induces by $\rho^{\sharp}$
	\[\bar{\rho}:C\left[(T_{r,s})_{1\leq r,s\leq n}, (\det)^{-1} \right]/I\simeq
		C\left[(Y_{i,j})_{1\leq i,j\leq m}, (\det)^{-1} \right]/J.
\]
Let $t_{r,s}$ be the image of $T_{r,s}$ in the quotient on the left, and recall that $y_{i,j}$ are that of 
$Y_{i,j}$ in the quotient on the right, then the preceding formula can be written
\[
	\bar{\rho}:C\left[(t_{r,s})_{1\leq r,s\leq n}\right]\simeq
	C\left[(y_{i,j})_{1\leq i,j\leq m}\right].
\]
Set ${\goth V}=\Spec \left(C\left[(t_{r,s})_{1\leq r,s\leq n}\right]\right)$, this is an algebraic subgroup of $\GL_n(C)$, it is isomorphic to ${\goth U}$ via the morphism induces by $\bar{\rho}$, denoted by abuse of language
$\rho :{\goth U}\simeq {\goth V}$. \par

The composed morphism (see (\ref{torsor}))
\begin{equation} \label{3}
	\ph :{\cal T}\otimes_{K}\widetilde{K}\stackrel{\psi^{-1}}{\longrightarrow}{\goth U}\times_{C}\widetilde{K}
	\stackrel{\rho \times \id_{\widetilde{K}}}{\longrightarrow} {\goth V}\times_{C}\widetilde{K}
\end{equation}
is an isomorphism of $\widetilde{K}$-schemes, equivariant for the actions of ${\goth U}$ and ${\goth V}$, this means that for any $\si $ in ${\goth U}$ we have 
$\ph \circ \si^{\flat}=\la_{\rho(\si )}\circ \ph$, where, as before, 
$\la_{\rho(\si )}$ is the endomorphism of ${\goth V}\times_{C}\widetilde{K}$ coming from the right translation by
$\rho(\si )$ on ${\goth V}$ (\ref{Galois action}).

\begin{lem}  Let $\ph^{\sharp}$ be the comorphism of $\ph $ (see (\ref{2}))and for any $r,s=1,\cdots ,n$ set
$z_{r,s}= \ph^{\sharp}(t_{r,s})$ (recall that ${\goth V}=\Spec C[(t_{r,s})_{1\leq r,s\leq n}]$). Then, for all $\si \in {\goth U}$, there exists a matrix $a(\si )\in \GL_n(C)$
such that we have the equality of matrices: $(\si(z_{r,s}))_{1\leq r,s\leq n}=(z_{r,s})_{1\leq r,s\leq n}a(\si )$.
\end{lem}
\begin{proof} Denote by $\la_{\rho (\si)}^{\sharp}$ the comorphism of the right translation by $\rho (\si)$ on
${\goth V}\times_{C}\widetilde{K}$, we have the equalities of matrices
	\[
\begin{array}{rl}
\left(\si(z_{r,s})\right)_{1\leq r,s\leq n} & =\left(\si \left(\ph^{\sharp}(t_{r,s})\right)\right)_{1\leq r,s\leq n}\\
	{} & =\left(\ph^{\sharp}\left(\la_{\rho (\si)}^{\sharp}(t_{r,s})\right)\right)_{1\leq r,s\leq n},
\end{array}
\]
because $\ph $ is equivariant, and
	\[\left(\la_{\rho (\si)}^{\sharp}(t_{r,s})\right)_{1\leq r,s\leq n}=(t_{r,s})_{1\leq r,s\leq n}a(\rho(\si ))
\]
where for any $\tau \in {\goth V}$ the matrix $a(\tau )$ is in $\GL_n(C)$ and is such that the formula 
$\left(\tau (t_{r,s})\right)_{1\leq r,s\leq n}=\left(t_{r,s}\right)_{1\leq r,s\leq n}a(\tau )$ defines the images of the $t_{r,s}$ by the comorphism $\la_{\tau}^{\sharp}$ of the right translation on ${\goth V}$ by $\tau $. We have find
	\[\left(\si(z_{r,s})\right)_{1\leq r,s\leq n}=\left(z_{r,s}\right)_{1\leq r,s\leq n}a(\rho(\si ))
\]
with $a(\rho(\si ))$ in $\GL_n(C)$.
\end{proof} 
The fact that $\ph $ is an isomorphism implies that $ R\otimes_{K}\widetilde{K}$ is generated over
$\widetilde{K}$ by the $z_{r,s}$, $1\leq r,s\leq n$, indeed $R\otimes_{K}\widetilde{K}$ is generated over
$\widetilde{K}$ by the $C$-space $V:=\sum_{1\leq r,s\leq n}Cz_{r,s}$ and the lemma shows that this space $V$ is (globally) invariant under the action of the Galois group ${\goth U}$. The (ordinary) Galois group 
$Gal(\widetilde{K} /K)$ acts as usual on the right hand factor of $R\otimes_{K}\widetilde{K}$ and trivially on the left one, then we see that $R$ is generated over $K$ by the $z_{r,s}$, $1\leq r,s\leq n$.\par

Another consequence of the previous lemma is that the matrix
\[D\stackrel{\mathrm{def}}{=}\left(z'_{r,s}\right)_{1\leq r,s\leq n}\left(z_{r,s}\right)^{-1}_{1\leq r,s\leq n} ,
\]
is in $\M_{n}(K)$, 
then, because $\ph^{\sharp}$ is an isomorphism, the ring $R$ is generated by the entries of a fundamental matrix of the equation $Y'=DY$, we know also that $R$ is a simple differential ring.  It follows that $R$, resp. $E$, is the Picard-Vessiot ring, resp. field, over $K$ of this equation.\par
 To a rational representation $\rho :dGAl(L/K)\rightarrow \GL_n(C)$ we have associated an element $[D]$ of 
 $Z_n(L/K)$, this is clearly the inverse map of $[A]\mapsto \rho_{[A]}$.
\end{proof}

\section{A correspondance.}
Let $K^{\mathrm{diff}}$ be a universal Picard-Vessiot extension of $K$ and set $G^{\mathrm{diff}}=dGal(K^{\mathrm{diff}}/K)$.We choose once of all an identification $\GL_n(C)=\GL(C^{n})$.\par

Let $\underline{\Rep}_n(G^{\mathrm{diff}})$ be the category of representations of $G^{\mathrm{diff}}$ in $\GL_n(C)$: the objects are morphisms $\rho :G^{\mathrm{diff}}\rightarrow \GL_n(C)$, an arrow $f: \rho_1\rightarrow \rho_2$ is a $C$-linear map from $C^{n}$ into itself such that, for any $g\in G^{\mathrm{diff}}$, the following diagram is commutative
	\[
\begin{array}{ccc}
	C^{n} & \stackrel{\rho_1(g)}{\longrightarrow} & C^{n} \\
	f\downarrow & {} & \downarrow f\\
	C^{n} & \stackrel{\rho_2(g)}{\longrightarrow}  & C^{n}
\end{array}
\] \par

To define the category $\underline{Z}_n(K)$ we need the following remarks. Let $M$ and $N$ be two elements of $\M_n(K)$, we say that they are equivalent if there exists $U$ and $V$ in $\GL_n(K)$ such that
$N=VMU$. We denote by $\overline{M}$ the equivalent class of $M$.
Let $A_i\in \M_n(K)$, $i=1,2$ and let $M\in \M_n(K)$ such that 
\begin{equation} \label{M'}
	M'=A_2M-MA_1.
\end{equation}
Let $B_i\in [A_i]$, let $U_i\in \GL_n(K)$ such that
	\[A_i=U'_iU^{-1}_i+U_iB_iU^{-1}_i,
\]
then an easy calculation shows that
	\[(U^{-1}_2MU_1)'=B_2(U^{-1}_2MU_1)-(U^{-1}_2MU_1)B_2.
\]
Suppose that $M\in \GL_n(K)$ and satisfies (\ref{M'}), then
	\[(M^{-1})'=A_1M^{-1}-M^{-1}A_2.
\]
Now we can define the category $\underline{Z}_n(K)$. Its objects are elements of $Z_n(K)$ (see (\ref{def Zn})), an arrow $[A_1]\rightarrow [A_2]$, where $A_1$ and $A_2$ are elements of $\M_n(K)$, is an equivalence class $\overline{M}$ in $M_n(K)$ such that there exists $M\in \overline{M}$ satisfying (\ref{M'}). The two preceding formulas show that this definition does not depend on the choice of $A_i$ in $[A_i]$, $i=1,2$, and that invertible arrows in $\underline{Z}_n(K)$ correspond to equivalence classes of invertible matrices. We explain the composition of arrows. Let $\overline{M}:[A_1]\rightarrow [A_2]$ and $\overline{N}:[A_2]\rightarrow [A_3]$ two arrows of 
$\underline{Z}_n(K)$, choose $M\in \overline{M}$, $N\in \overline{N}$ such that
	\[M'=A_2M-MA_1 \ \ \text{ and } \ \ N'=A_3N-NA_2,
\]
then we see that 
	\[(NM)'=A_3NM-NMA_1.
\]
The composed arrow is $\overline{N}\circ \overline{M}=\overline{NM}$, for a good choice of representing elements
of the different classes of matrices.
\par
Then $\underline{Z}_n(K)$ is a category, indeed it is easily to see that it is an additive category.

\begin{thm} \label{2main thm}
The two categories $\underline{Z}_n(K)$ and $\underline{\Rep}_n(G^{\mathrm{diff}})$ are equivalent. On objects, this equivalence is $[A]\mapsto c_{[A]}$ (see (\ref{c_A})).
\end{thm}
\begin{proof} Note that here to write $c_{[A]}$ is an abuse of notation, if $L/K$ is the Picard-Vessiot extension 
(contained in $K^{\mathrm{diff}}$) associated to the equation $Y'=AY$, we denote always $c_{[A]}$ the representation
	\[G^{\mathrm{diff}}\stackrel{\mathrm{restriction}}{\longrightarrow} dGal(L/K)
	\stackrel{c_{[A]}}{\longrightarrow} \GL_n(C).
\]
The map $[A]\mapsto c_{[A]}$ on objects of the categories has been constructed in the previous theorem, it is one to one. Let $[A_1]$ and $[A_2]$ be two objects of $\underline{Z}_n(K)$ and $\overline{M}:[A_1]\rightarrow [A_1]$ be an arrow, select $M\in \overline{M}$ such that $M'=A_2M-MA_1$. Let $F_1,F_2\in \GL_n(K^{\mathrm{diff}})$ be fundamental matrices for respectively the equations $Y'=A_1Y$ and $Y'=A_2Y$.  Then $F'_i=A_iF_i$, $i=1,2$.
Let $f=F^{-1}_2MF_1$, a priori $f$ is in 
$\GL_n(K)$, but
	\[
\begin{array}{rcl}
	f' & =(F^{-1}_2)'MF_1+F^{-1}_2M'F_1+F^{-1}_2MF'_1 \\
	{} & =(-F^{-1}_2A_2)MF_1+F^{-1}_2(A_2M-MA_1)F_1\\
	{} & +F^{-1}_2MA_1F_1 =0.
\end{array}
\]
Then $f=F^{-1}_2MF_1$ is in $\GL_n(C)$. Now we prove that $f$ is a morphism from $c_{[A_1]}$ to $c_{[A_2]}$. Let 
$g$ be an element of $G^{\mathrm{diff}}$. Applying $g$ to the relation $f=F^{-1}_2MF_1$ we find
	\[f=g(F^{-1}_2)Mg(F_1)=g(F^{-1}_2)F_2fF^{-1}_1g(F_1)=c_{[A_2]}(g)^{-1}fc_{[A_1]}(g),
\]
(see (\ref{c_A})) for all $g$. This means that $f:c_{[A_1]}\rightarrow c_{[A_2]}$ is a map in $\underline{\Rep}_n(G^{\mathrm{diff}})$.\par

Conversely let $f:\rho_1\rightarrow \rho_2$ be an arrow of $\underline{\Rep}_n(G^{\mathrm{diff}})$, then we can see $f$ as a matrix with coefficient in $C$. 
We know that there exists $A_i$ in $\M_n(K)$ such that $\rho_i=c_{[A_i]}$, $i=1,2$. Let as before $F_i$ be a fundamental matrix for the equation $Y'=A_iY$. Set $M=F_2fF^{-1}_1$. \par
\noindent - We prove that $M$ is in $\M_n(K)$.
 The fact that $f$ is a morphism of representations means that for all $g$ in $G^{\mathrm{diff}}$ we have
	\[fc_{[A_1]}(g)=c_{[A_1]}(g)f,
\]
which is equivalent to
	\[fF^{-1}_1g(F_1)=F^{-1}_2g(F_2)f,
\]
then 
	\[F_2fF^{-1}_1=g(F_2)fg(F^{-1}_1)=g(F_2fF^{-1}_1).
\]
This prove that the entries of $M$ are in $K$.\par
\noindent - We prove the formula $M'=A_2M-MA_1$. We have
	\[M'=F'_2fF^{-1}_1+F_2f(F^{-1}_1)'=A_2F_2fF^{-1}_1+F_2f(-F^{-1}_1A_1)
\]
which is the expected formula.
\end{proof}


\end{document}